\documentclass[12pt,a4paper]{article}
\usepackage{a4wide,amsmath,amsthm,amssymb,latexsym}

\newtheorem{theorem}{Theorem}
\newtheorem{lemma}[theorem]{Lemma}
\newtheorem{corollary}[theorem]{Corollary}
\newtheorem{proposition}[theorem]{Proposition}

\newcommand{\N}{{\mathbb N}}
\newcommand{\Z}{{\mathbb Z}}
\newcommand{\F}{{\mathbb F}}
\newcommand{\FF}{{\mathcal F}}

\newcommand{\rank}{\mathop{rank}}

\newcommand{\imm}{\mathrm{Im}}

\renewcommand{\dim}{\mathrm{dim}}

\newcommand{\GL}{\operatorname{GL}}
\newcommand{\sr}{\operatorname{SupR}}
\newcommand{\Ann}{\operatorname{Ann}}

\renewcommand{\phi}{\varphi}

\newcommand{\sig}{\operatorname{sig}}
\newcommand{\nd}{\operatorname{nd}}

\begin{document}

\title{On the semigroup of square matrices}
\author{Ganna Kudryavtseva and Volodymyr Mazorchuk}
\date{}
\maketitle

\begin{abstract}
We study the structure of nilpotent subsemigroups in the semigroup
$M(n,\F)$ of all $n\times n$  matrices over a field, $\F$, with
respect to the operation of the usual matrix multiplication. We
describe the maximal subsemigroups among the nilpotent subsemigroups
of a fixed nilpotency degree and classify them up to isomorphism.
We also describe isolated and completely isolated subsemigroups and
conjugated elements in $M(n,\F)$.
\end{abstract}

\section{Introduction}\label{s1}

The structure and combinatorics of the classical finite transformation
semigroups is now relatively well understood. A lot of information about
such semigroups as the full finite inverse symmetric semigroup
$\mathcal{I}\mathcal{S}_n$, the full transformation semigroup 
$\mathcal{T}_n$ and the semigroup $\mathcal{P}\mathcal{T}_n$ of all 
partial transformations on $\{1,2,\dots,n\}$ can be found for example
in the monographs \cite{CP,Hi,Ho,La,Li} or other numerous papers studying 
transformation semigroups. Two directions of study for 
transformation semigroups, which developed over the last 15 years, 
is the study of conjugated elements and nilpotent subsemigroups in 
these semigroups, which  resulted into several nice structural 
and combinatorial results, see \cite{GK0,GK1,GK2,GK3,GTS,GM,KM}.

Certainly, the passage to infinite transformation semigroups completely 
changes the picture. However, it is still possible to obtain some information
under reasonable ``finiteness'' conditions. One of the most classical examples 
of an infinite object possessing several properties, inherent in finite 
objects, is the algebra of all linear operators on a finite-dimensional 
vector space over a field.  Forgetting the addition of linear operators
one gets a semigroup, isomorphic to the semigroup $M(n,\F)$ of all 
$n\times n$ square matrices with coefficients from a field $\F$ with respect
to the operation of usual matrix multiplication. This semigroup has also
been studied by many authors, but so far not as intensively as the classical
transformation semigroups. Many interesting results about $M(n,\F)$ can be
found in the recent monograph \cite{Ok2} and in its references.

The aim of the present paper is to contribute to the study of  $M(n,\F)$ with
results in several directions. The main emphasis is made on the study of
conjugated elements and nilpotent subsemigroups of $M(n,\F)$. In the appendix 
we also address the problem of the study of isolated and completely 
isolated subsemigroups. These different directions are not immediately 
related with each other, however, when combined together, they give a very 
nice illustration for the fact that in the passage from the classical
transformation semigroups to matrix semigroups one should expect that 
some results would be transfered to very similar results, while some other
results would sound quite differently. And even in the case of similar 
results, the passage to matrix semigroups substantially increases the level 
of technical difficulties.

The notion of conjugated elements in group theory is very important and has
a lot of applications (for example for the study of automorphisms, characters
or representations). There are several ways to extend this notion for
semigroups. Two most straightforward generalizations are: conjugation with
respect to an invertible element, and the transitive closure of the 
$ab\sim ba$ relation. Both these notions provide some invariant on semigroups
and hence can be applied for the study of automorphisms, endomorphisms and
representation. The already existing literature, where these notions were
studied (see for example \cite{GK0,Li,KM}) suggests that the second 
generalization is more interesting than the first one. For $M(n,\F)$ the
description of equivalence classes with respect to the conjugation by
invertible matrices is a classical problem, the answer to which is given by
the Jordan normal form of a matrix (in the case of an algebraically closed
field). In Section~\ref{sn6} of the present paper we describe the  equivalence
classes with respect to the transitive closure of the  $ab\sim ba$ relation.
It turns out that these are given by the ``invertible part'' of the Jordan
normal form. This is very similar to the results obtained in \cite{GK0,KM}
for transformation semigroups. The results of Section~\ref{sn6} of the 
present paper and the results of \cite{GK0,KM} were the main motivation for 
the abstract approach to the study of conjugation for semigroups, developed
in \cite{Ku}.

The majority of the paper is devoted to the study of nilpotent subsemigroups
of $M(n,\F)$. To start with, we would like to remark that the notion of  a
nilpotent semigroup has been used in the literature in at least  three
different senses. The most classical one is the notion of {\em nilpotent
semigroup in the sense of Maltsev}, \cite{Ma}, which is defined by means of 
the identities for nilpotent groups, rewritten without $g^{-1}$ terms.  The
matrix semigroups, nilpotent in the sense of Maltsev, were recently studied 
in \cite{O}. However, in this paper we are going to use another notion of a
nilpotent semigroup, which comes from the ring theory. A semigroup, $S$,  with
the zero element $0$ is called {\em nilpotent} of {\em nilpotency degree}  
$\nd(S)=k$ provided that $a_1a_2\dots a_k=0$ for any $a_1,\dots,a_k\in S$ while
there exist $b_1,\dots,b_{k-1}\in S$ such that $b_1\dots b_{k-1}\neq 0$. This
notion is almost as old as the first one and goes back at least till \cite{Shev}.  From now on we will use only the last notion of nilpotent
semigroups.

The study of nilpotent subsemigroups of certain semigroups of partial
transformations, in particular, of the semigroup ${\mathcal IS}_n$, was
originated in \cite{GK3}. It happened that the combinatorial data, describing
the maximal nilpotent subsemigroups of finite transformation semigroups, is
usually a certain partial order on the underlined set, on which the semigroup
acts. This philosophy was successfully used later in \cite{GK1,GK2} and
generalized on the infinite case in \cite{Sh3,Sh1}. In \cite{GTS} the
combinatorial description of maximal nilpotent subsemigroups of 
${\mathcal IS}(M)$ was used to determine the group of automorphisms for 
these subsemigroups. 

In the present paper we generalize this technique to study nilpotent
subsemigroups of the semigroup $M(n,\F)$. A part of this study is very easy.
It is well-known that the maximal nilpotent subalgebras of $M(n,\F)$ 
correspond bijectively to complete flags in $\F^n$. An elementary exercise
in linear algebra shows that a maximal nilpotent subsemigroup of $M(n,\F)$
is in fact a subalgebra. Hence maximal nilpotent subsemigroup of $M(n,\F)$
are also described by complete flags in $\F^n$.
This is a perfect analogy with the combinatorial data describing
maximal nilpotent subsemigroups in transformation semigroups. So, there is no
problem to classify such subsemigroups, what we do in Section~\ref{s2}, where
we also collect some other statement necessary for or main goal: to classify
maximal nilpotent subsemigroups of $M(n,\F)$ up to isomorphism. At this point
we face our main difficulty in comparison to finite transformation subgroups.
The fact that the cardinality of $M(n,\F)$ is infinite for infinite $\F$ makes
it impossible to distinguish non-isomorphic maximal nilpotent subsemigroups of
$M(n,\F)$ by their cardinality, as it was done in case of ${\mathcal IS}_n$ 
in \cite{GK2}. Therefore we classify all maximal nilpotent subsemigroups of
$M(n,\F)$ up to isomorphism in Section~\ref{s3} using completely different
arguments.

In the Appendix we classify all isolated and completely isolated 
subsemigroups of  $M(n,\F)$ for finite $\F$. This is not directly related 
to the main content  of the paper, however, from our point of view, these
results, especially for isolated subsemigroups, are not immedeately 
expected, are non-trivial and  again give a very nice (and not yet
well-understood) illustration of  similarity with the classical transformation
semigroups. The results of the Appendix were also a part of the  motivation  
for the study of isolated subsemigroups in the variant of $\mathcal{T}_n$, 
see \cite{MT}.


\section{Conjugated elements}\label{sn6}

Let $S$ be a monoid and $G$ be its group of units. 
The elements $x, y\in S$ are said to be
$G$-conjugated provided there exists $g\in G$ such that
$x=g^{-1}yg$. This is denoted by $x\sim_G y$. The binary relation
$\sim_G$ is an equivalence relation on $S$.
The elements $a$ and $b$ of a semigroup, $S$, are called
{\em primarily $S$-conjugated} if there exist such $x, y\in S$ that $a =xy$
and $b=yx$. The binary relation $\sim_{pS}$ of primary  $S$-conjugation
is reflexive and symmetric, but not transitive in general. We denote
by $\sim_S$ the transitive closure of this relation. If $x\sim_S y$,
the elements $x$ and $y$ will be called {\em $S$-conjugated.}
Both $\sim_G$ and $\sim_S$ generalize the notion of
conjugated elements in a group, and in the general case the relations
$\sim_G$ and $\sim_S$ do not coincide. However, $\sim_G$ is always a
subset of $\sim_S$.

The group of units in $M(n, \F)$ is $\GL(n, \F)$.
If $\F$ is algebraically closed, then the description of
$\GL(n, \F)$-conjugated elements in $M(n, \F)$ is
the classical Jordan theorem of the basic linear algebra:
$A\sim_{\GL(n, \F)} B$ if and only if the Jordan normal forms of $A$
and $B$ coincide (up to a permutation of Jordan cells).

Let $A\in M(n,\F)$. For each $k\geq 0$ we have the inclusion
$\F^n\supset A^k(\F^n)\supset A^{k+1}(\F^n)\supset 0$, which implies
$n\geq\dim(A(\F^n))\geq\dots\geq\dim(A^k(\F^n))\geq\dots\geq 0$.
Since only at most $n$ of these inequalities can be strict, we can
assert that, starting from some power $t$, we have
$A^t(\F^n)=A^{t+i}(\F^n)$
for each $i\geq 0$. Denote by $A_s$ the linear operator defined as
follows: $A_s(v)=A(v)$, $v\in A^t(\F^n)$; $\ker(A_s)=\ker(A^t)$.
Remark that for an algebraically closed field, $\F$, if we fix a
Jordan basis for $A$, then $A_s$ is obtained from $A$ by replacing
all Jordan cells with eigenvalue $0$ by zero blocks. The main result
of this section is the following.

\begin{theorem}\label{maincon}
Two matrices $A$ and $B$ from $M(n,\F)$ are
$M(n,\F)$-conjugated if and only if $A_s$ and $B_s$ are
$\GL(n,\F)$-conjugated.
\end{theorem}

\begin{proof}{\em Sufficiency.} We start with the following lemma:

\begin{lemma}\label{conl2}
$A\sim_{M(n, \F)} A_s$ for every $A\in M(n, \F)$. In particular,
$A\sim_{M(n,\F)}0$ for every nilpotent $A\in M(n,\F)$.
\end{lemma}

\begin{proof}
Denote by  $V_i=\imm(A^i)$, $i\in\N$, $V_0=\F^n$. Clearly, $V_i\subset
V_{i-1}$ for all $i$. Let $V'_i$ be a complement to $V_i$, $i\in\Z_+$,
such that $V'_i=\ker(A^i)$ for all $i\geq t$.
Set $B_i=e(V_i,V'_i)Ae(V_{i-1},V'_{i-1})$. Then $B_1=A$ and
$B_n=A_s$. Hence, to complete the proof it is enough to show
that $B_i\sim_{pM(n, \F)}B_{i+1}$ for all $i$. For a fixed $i$ denote
$u_i=e(V_i,V'_i)$ and $v_i=e(V_i,V'_i)Ae(V_{i-1},V'_{i-1})=B_i$.
Then $u_iv_i=v_i=B_i$. Moreover,
\begin{multline*}
v_iu_i=e(V_i,V'_i)Ae(V_{i-1},V'_{i-1})e(V_i,V'_i)=\\=
e(V_i,V'_i)Ae(V_i,V'_i)=e(V_{i+1},V'_{i+1})Ae(V_i,V'_i)=B_{i+1}.
\end{multline*}
\end{proof}

Let $A_s$ and $B_s$ be $\GL(n, \F)$-conjugated. By Lemma~\ref{conl2}, we have  
$A\sim_{M(n, \F)} A_s$ and $B\sim_{M(n, \F)} B_s$. Hence the transitivity of 
the relation $\sim_{M(n, \F)}$ implies $A\sim_{M(n, \F)} B$.

{\em Necessity.}  Clearly, it is enough to  consider only the case
when  $A$ and $B$ are  primarily $M(n, \F)$-conjugated. So, let
$A=XY$ and  $B=YX$ for some $X, Y\in M(n, \F)$. Let
$V_1=\imm(A^n)$, $W_1=\imm(B^n)$, $V_2=\ker(A^n)$ and
$W_2=\ker(B^n)$. Then we have $\F^n=V_1\oplus V_2=W_1\oplus W_2$
and
\begin{equation}\label{eqeqnew}
\rank(A^n)=\rank(A^{i})\quad\text{ and }\quad
\rank(B^n)=\rank(B^{i})\quad  \text{ for all }\quad  i>n.
\end{equation}
Since $B(\ker(X))=0$ and $B(W_1)=W_1$ it follows that $\ker(X)\cap
W_1=0$, analogously $\ker(Y)\cap V_1=0$. Let $x\in V_2$. Then
$A^{n+1}x=XB^nYx=0$, and $\ker(X)\cap W_1=0$ forces $Yx\in W_2$,
that is $Y:V_2\to W_2$. Analogously one shows that $X:W_2\to V_2$.
If $x\in V_1$ then $x=A^nx'$ for some $x'\in V_1$ since
$A^n:V_1\to V_1$ is bijective by \eqref{eqeqnew}. Hence
$Yx=YA^nx'=B^nYx'\in W_1$. This implies $Y:V_1\to W_1$ and,
analogously,  $X:W_1\to V_1$. The equalities $\ker(X)\cap W_1=0$
and $\ker(Y)\cap V_1=0$ now say that $Y:V_1\to W_1$ and $X:W_1\to
V_1$ are bijections. In particular,
\begin{equation}\label{eqeqnew2}
Ye=fYe\quad \text{ and }\quad Xf=eXf,
\end{equation}
where $e=e(V_1,V_2)$ and $f=e(W_1,W_2)$.

Further, we have $A_s=eAe$ as the actions of these operators on both
$V_1$ and $V_2$ coincide. Analogously $B_s=fBf$. Now from
$A_s=eAe=eXYe$, $B_s=fBf=fYXf$ and \eqref{eqeqnew2}
we derive $A_s=eXYe=eXfYe=eXffYe$
and $B_s=fYXf=fYeXf=fYeeXf$, implying, in particular, that $A_s$ and
$B_s$ are primarily $M(n, \F)$-conjugated.

The fact that $Y$ induces a bijection from $V_1$ to $W_1$ implies that
$\dim(V_1)=\dim(W_1)$ and thus $\dim(V_2)=\dim(W_2)$. Let
$Z$ be any matrix such that $\ker(Z)=V_1$ and $Z:V_2\to W_2$ is a
bijection, which exists since $\dim(V_2)=\dim(W_2)$. Consider the
matrix $M=fYe+Z$. We have $M(V_1)=fYe(V_1)=W_1$ since
$Z(V_1)=0$ and $M(V_2)=Z(V_2)=W_2$ since $e(V_2)=0$. Hence
$fYe+Z\in \GL(n,\F)$.

From the definition of $Z$ we obtain $fZ=Ze=0$, which implies
\begin{displaymath}
B_sM= fYeeXf(fYe+Z)=fYeeXffYe=(fYe+Z)eXffYe=MA_s.
\end{displaymath}
Hence $B_s=MA_sM^{-1}$, which completes the proof.
\end{proof}

We would like to remark that, for an algebraically closed $\F$,
Theorem~\ref{maincon} provides a criterion of
$M(n,\F)$-conjugacy in terms of Jordan normal forms:
$A$ and $B$ are $M(n,\F)$-conjugated if and only if their
Jordan normal forms coincide up to a permutation of Jordan cells
and a replacement of some nilpotent Jordan cells by some
other nilpotent Jordan cells.


\section{Nilpotent subsemigroups of $M(n,\F)$}\label{s2}

The results of this section are mostly known or can be easily deduced
from the existing literature. However, we did not manage to single them
out in the literature in one piece and in the form we present them. 
Our presentation emphasizes the similarity with the corresponding 
results for transformation semigroups, see \cite{GK1,GK2,GK3,GM}.
We include the proofs for the sake of completeness.

\begin{lemma}\label{l1}
Let $T$ be a maximal nilpotent subsemigroup of $M(n,\F)$.
Then $T$ is, in fact, a subalgebra of $M(n,\F)$.
\end{lemma}

\begin{proof}
Let $T$ be a maximal nilpotent subsemigroup of $M(n,\F)$.
Denote by  $T'$ the linear span of $T$ inside $M(n,\F)$.
Then every element in $T'$ is a finite
linear combination of elements from $T$, and $T\subset T'$. Since
matrix multiplication is bilinear and $T$ is a semigroup, we
obtain that $T'$ is in fact a subsemigroup of $M(n,\F)$. To complete
the proof it is enough to show that $T'$ is nilpotent, which
follows immediately from the nilpotency of $T$ and
\begin{displaymath}
\left(\sum_{{i_1}=1}^{m_1}\lambda_{i_1}^{(1)} a_{i_1}^{(1)}\right)\dots
\left(\sum_{{i_s}=1}^{m_s}\lambda_{i_s}^{(s)} a_{i_s}^{(s)}\right)=
\sum_{{i_1}=1}^{m_1}\dots\sum_{{i_s}=1}^{m_s}\left(
\lambda_{i_1}^{(1)}\dots \lambda_{i_s}^{(s)}
a_{i_1}^{(1)}\dots a_{i_s}^{(s)}\right).
\end{displaymath}
\end{proof}

Nilpotent subalgebras of $M(n,\F)$ are classical objects both in algebraic
geometry and in representation theory and have a natural description in terms 
of flags in $\F^n$. By a {\em flag}, $\FF$, in $\F^n$ of {\em length}
$l(\FF)=k$ we will mean a filtration, $\FF:0=V_0\subset V_1\subset\dots 
\subset V_k=\F^n$, of $\F^n$ by subspaces, such that $V_i\neq V_{i+1}$, $i=1,\dots,k-1$. The tuple $(d_1,\dots,d_k)$, where $d_i=\dim(V_i/V_{i-1})$,
will be called the {\em signature} of $\FF$ and will be denoted by $\sig(\FF)$.
A basis, $\{v_1,\dots,v_n\}$, of $\F^n$ is said to be an {\em $\FF$-basis} 
if $v_1,\dots,v_{\dim(V_i)}$ is a basis of $V_i$ for all  $i=1,\dots, k$.

If $S$ is a semigroup with $0$ and $k\geq 2$, we denote by $N_k(S)$ the set of
all nilpotent subsemigroups of $S$ of nilpotency degree $k$. This is a 
partially ordered set with respect to inclusion and the set $\cup_{k\geq 2} 
N_k(S)$ is the set of all (non-zero) nilpotent subsemigroups of $S$, which is
also partially ordered with respect to inclusion.

To each flag $\FF$ in $\F^n$ of length $k$, we associate the nilpotent
subsemigroup $\phi(\FF)=\{a\in M(n,\F)|a(V_i)\subset V_{i-1}\text{ for 
all }i\}$ of $M(n,\F)$ of nilpotency degree $k$. To each nilpotent
subsemigroup $S$ of $M(n,\F)$ of nilpotency degree $k$, we associate
the flag $\psi(S)$ in $\F^n$ defined as follows: $0\subset \langle
S^{k-1}(\F^n)\rangle\subset \langle S^{k-2}(\F^n)\rangle\subset \dots
\subset \langle S(\F^n)\rangle\subset \F^n$. The correctness of
this definition follows from the following statement.

\begin{proposition}\label{l2}
\begin{enumerate}
\item If $\FF$ is a flag in $\F^n$ of length $k$, then $\phi(\FF)$ is a
nilpotent semigroup of $M(n,\F)$ of nilpotency degree $k$.
\item If $S$ is a nilpotent subsemigroup of $M(n,\F)$ of nilpotency
degree
$k$, then $\psi(S)$ is a flag of $F^n$.
\end{enumerate}
\end{proposition}

\begin{proof}
Obviously, $\phi(\FF)$ is a nilpotent semigroup and $\phi(\FF)^k=0$.
Let $B=\{v_1,\dots,v_n\}$ be some $\FF$-basis in $\F^n$.
Consider the element $a\in M(n,\F)$, defined on $B$ as follows:
$a(v_{\dim(V_{i})})=v_{\dim(V_{i-1})}$, $i=2,\dots,k$, and $a(v_j)=0$
otherwise. Then $a(V_i)\subset V_{i-1}$ by the definition and hence
$a\in\phi(\FF)$. But $a^{k-1}(v_n)= v_{\dim(V_1)}$ and hence
$a^{k-1}\neq 0$. This means that $\nd(\phi(\FF))>k-1$
and thus equals $k$. This proves the first statement.

Since $a_1a_2a_3\dots a_{i+1}=(a_1a_2)a_3\dots a_{i+1}$ we have
$S^{i+1}\subset S^i$ and thus obtain the inclusion
$\langle S^{i+1}(\F^n)\rangle\subset \langle S^i(\F^n)\rangle$. Hence
it is only left to prove that $\langle S^i(\F^n)\rangle\neq \langle
S^{i+1}(\F^n)\rangle$. But $\langle S^{i+1}(\F^n)\rangle$ has a basis,
consisting of elements from $S^{i+1}(\F^n)$. Further,
$S^k=0$ implies $S^{k-i-1}(S^{i+1}(\F^n))=0$ and, using the fact that
the elements from $S$ are linear operators on $\F^n$, we obtain that
$S^{k-i-1}(\langle S^{i+1}(\F^n)\rangle)=0$. To complete the proof it
is now enough to show that $S^{k-i-1}(\langle S^{i}(\F^n)\rangle)\neq
0$. Since $S^{i}(\F^n)\subset \langle S^{i}(\F^n)\rangle$,
it is even enough to show that
$S^{k-i-1}(S^{i}(\F^n))\neq 0$, which reduces to $S^{k-1}(\F^n)\neq
0$. The last inequality follows from the fact that $\nd(S)=k$ and
therefore
$S^{k-1}\neq 0$.
\end{proof}

Proposition~\ref{l2} immediately implies the following classical results:

\begin{corollary}\label{cnew1}
If $S$ is a nilpotent subsemigroup of $M(n,\F)$, then $\nd(S)\leq n$.
\end{corollary}

\begin{corollary}\label{cnewnew2}
A subsemigroup, $S\subset M(n,\F)$, is nilpotent if and only if $S$
consists of nilpotent elements.
\end{corollary}

\begin{lemma}\label{llnew}
\begin{enumerate}
\item Let $i\in\{2,\dots,k\}$. Then for every $v\in V_i\setminus
V_{i-1}$ and $w\in V_{i-1}\setminus V_{i-2}$ there exists
$a\in \phi(\FF)$ such that $a(v)=w$.
\item $\phi(\FF)^i(\F^n)= V_{k-i}$ for all $i=1,\dots,k$.
\item $\phi(\FF)$ is a subalgebra of $M(n,\F)$.
\end{enumerate}
\end{lemma}

\begin{proof}
We can certainly find an $\FF$-basis $B$ of $V$ such that
both $v$ and $w$ belong to $B$. The linear operator $a\in M(n,\F)$,
which sends $v$ to $w$ and annihilates all other elements of $B$,
belongs to $\phi(\FF)$, which proves the first statement.
The second and the third statements are obvious.
\end{proof}

Let $\FF$ be a flag in $\F^n$. We will say that $\sig{\FF}$ is also
the {\em signature} of the semigroup $\phi(\FF)$.
Let $k$ be a positive integer, $k\geq 2$. A nilpotent subsemigroup,
$S\subset M(n,\F)$, is called {\em $k$-maximal} provided that $S$
is a maximal element in $N_k(S)$.

\begin{proposition}\label{pro3}
$\phi(\FF)$ is $k$-maximal with $k=l(\FF)$ and every $k$-maximal
nilpotent subsemigroup of $M(n,\F)$ equals $\phi(\FF)$ for some
flag, $\FF$, of $\F^n$ of length $k$.
\end{proposition}

\begin{proof}
Let $S$ be a nilpotent semigroup of $M(n,\F)$ of nilpotency degree $k$.
From the definitions of $\phi$ and $\psi$ it follows that $S\subset
\phi(\psi(S))$. Hence any $k$-maximal nilpotent subsemigroup of
$M(n,\F)$ has to be of the form $\phi(\FF)$ for a flag, $\FF$, of
$\F^n$ of length $k$.

Now assume that $S=\phi(\FF)\subset T$, where $T\in N_k(M(n,\F))$. Then
$S^i(\F^n)=V_{k-i}\subset T^i(\F^n)$ for all $i$ and hence it is enough to
prove that $S^i(\F^n)= T^i(\F^n)$. Assume that $j$ is minimal
such that this equality fails. This means that there exists $v\in
S^{j-1}(\F^n)=T^{j-1}(\F^n)$
and $a\in T$ such that $a(v)\not\in S^j(\F^n)$.
Consider a new flag, $\FF'$, which has the form
\begin{displaymath}
0=V_0\subset V_1\subset\dots\subset V_{k-j} \subset \langle T^j(\F^k)
\rangle\subset V_{k-j+1}\subset\dots V_k=\F^n.
\end{displaymath}
For every $i=0,1,\dots,k-1$,
$i\neq j-1$, pick $v_j\in S^i(\F^n)\setminus S^{i+1}(\F^n)$. Then the set
\begin{displaymath}
\{v_{k-1},\dots,v_{j},a(v),v,v_{j-2},\dots,v_0\}
\end{displaymath}
is linearly
independent,
as its components belong to the different subquotients of $\FF'$, so we
can
extend this set to a basis, say $B$, of $\F^n$, coordinated with
$\FF'$.
Now consider the elements $a_i$, $i\neq j-2,j-1$, defined as follows:
$a_i(v_i)=v_{i+1}$ and $a_i(w)=0$, $w\in B\setminus\{v_i\}$; and
elements $b,c$, defined as follows: $b(v_{j-2})=v$, $b(w)=0$, $w\in
B\setminus \{v_{j-2}\}$; $c(a(v))=v_{j}$, $c(w)=0$, $w\in B\setminus
\{a(v)\}$. It follows immediately that $b,c$ and all $a_i$ are
elements of $\phi(\FF)$ and hence of $T$.
But $a_{k-2}\dots a_{j}caba_{j-3}\dots a_1(v_1)=
v_{k-1}\neq 0$ and hence $\nd(T)>k$, a contradiction.
\end{proof}

\begin{corollary}\label{c4}
The maps $\phi$ and $\psi$ are mutually inverse bijections between
the set of $k$-maximal nilpotent subsemigroups of $M(n,\F)$
and the set of all flags in $\F^n$ of length $k$.
\end{corollary}

We will say that a flag $\FF$ is a {\em consolidation} of $\FF'$
provided that each $V'_i$ equals $V_j$ for some $j$. In particular,
if $\FF$ is a consolidation of $\FF'$ then either $\FF=\FF'$ or
$l(\FF)>l(\FF')$.

\begin{proposition}\label{pro5}
$\phi(\FF')\subset \phi(\FF)$ if and only if $\FF$ is a consolidation
of $\FF'$.
\end{proposition}

\begin{proof}
The ``if'' part follows from the definition of $\phi$. To prove the
``only if'' part we set $S=\phi(\FF')$, $T=\phi(\FF)$ and
actually have to show that for every $i$ there is $j$ such that
$S^i(\F^n)=T^j(\F^n)$. Assume that this is not the case and let
$i$ be minimal with this property.
From $S\subset T$ we have $S^i\subset T^i$.
Let $j$ be maximal such that $S^i(\F^n)\subset T^{j}(\F^n)$.
In particular, $S^i(\F^n)\neq T^{j}(\F^n)$ and $S^i(\F^n)\setminus
T^{j+1}(\F^n)\neq \varnothing$. Note that our choice of $i$ and $j$
guarantees that $T^{j}(\F^n)\setminus S^i(\F^n)\subset
S^{i-1}(\F^n)\setminus S^i(\F^n)$. Take
$v\in T^{j}(\F^n)\setminus S^i(\F^n)$. By the choice of $i$ and $j$,
$T^{j+1}(\F^n)\cap S^i(\F^n)$ is a proper subspace of $S^i(\F^n)$.
Since $S^{i+1}(\F^n)$ is a proper subspace of $S^i(\F^n)$ as well,
we get that $\left(T^{j+1}(\F^n)\cap S^i(\F^n)\right)\cup
S^{i+1}(\F^n)\neq S^i(\F^n)$ (a non-zero vector space is never
a union of two proper subspaces). Take
$w\in S^i(\F^n)\setminus \left(\left(T^{j+1}(\F^n)\cap
S^i(\F^n)\right)\cup S^{i+1}(\F^n)\right)$. By the first statement
of Lemma~\ref{llnew}, there exists $a\in S\subset T$ such that
$a(v)=w$, but $a(T^{j}(\F^n))\not\subset T^{j+1}(\F^n)$, as $v\in
T^{j}(\F^n)$ and $w\not\in T^{j+1}(\F^n)$ by the construction. This
gives a contradiction and completes the proof.
\end{proof}

\section{Isomorphism of maximal nilpotent subsemigroups}\label{s3}

In the present section we give a description of the $k$-maximal
nilpotent subsemigroups of $M(n,\F)$ for a fixed $k$ up to
isomorphism. 

\begin{corollary}\label{c7}
Let $\F$ be infinite. Then every two maximal $2$-nilpotent
subsemigroups of $M(n,\F)$ are isomorphic.
\end{corollary}

\begin{proof}
Clearly, $|M(n,\F)|=|\F|$ in case of infinite $\F$. Then, if
$a\in\phi(\FF)$ we get that $\lambda a\in \phi(\FF)$ for all
$\lambda\in\F$ and hence $|\F|\leq |\phi(\FF)|\leq |M(n,\F)|= |\F|$.
Thus $|\phi(\FF)|=|\F|$ and the statement follows. Now the
statement follows from the fact that two semigroups, $S$ and $T$, 
with zero multiplication (i.e $ab=0$ for all $a,b$) are isomorphic if 
and only if $|S|=|T|$ since any bijection $S\setminus\{0\}\to 
T\setminus\{0\}$ extends to an isomorphism via $0\mapsto 0$.
\end{proof}

It is clear that the last statement is completely wrong if we replace
``semigroup'' with ``algebra''. This shows the difference between these
two theories. The main reason for this difference is that in the case
of the ``algebra'' structure there is a nice finite invariant:
dimension. This invariant does not work in the ``semigroup''
case because vector spaces of different finite dimensions over an
infinite field still have the same cardinality.

\begin{lemma}\label{c8}
Let $\F=\F_q$ be a finite field, $|\F_q|=q=p^l$, where $p$ is a
prime. Let $\FF$ and $\FF'$ be two flags in $\F^n$ of length $2$.
The semigroups $S=\phi(\FF)$ and $T=\phi(\FF')$ are isomorphic
if and only if $\{\dim(V_1),\dim(V_2/V_1)\}=$
$\{\dim(V'_1),\dim(V'_2/V'_1)\}$.
\end{lemma}

\begin{proof}
As $\F$ is finite, all semigroups we consider will be finite as well,
so we can try to calculate their cardinalities. So, we assume
$\dim(V_1)=m$ and have  $\dim(V_2/V_1)=n-m$. Let $B=\{v_1,\dots,v_n\}$
be an $\FF$-basis of $\F^n$. An element, $a\in M(n,\F)$,
belongs to $S$ if and only if $a(v_i)=0$, $i=1,\dots,m$, and
$a(v_i)=\sum_{j=1}^m a_j v_j$, $a_j\in\F$, $i>m$. Hence
$|S|=q^{m(n-m)}$.

If $\dim(V'_1)=m'$, we get $|T|=q^{m'(n-m')}$. Further $|S|=|T|$
implies $q^{m(n-m)}=q^{m'(n-m')}$, which is equivalent to
$m(n-m)=m'(n-m')$. The left hand side of the last equality is
a quadratic polynomial in $m$, and the right hand side represents the
value of the same polynomial in the point $m'$. This implies that
there the equality holds if and only if $m=m'$ or $m=(n-m')$.
Now everything again follows from the fact that two semigroups, $S$ and 
$T$, with zero multiplication are isomorphic if and only if $|S|=|T|$.
\end{proof}

We switch to the case of arbitrary $\F$. The main result of this
section
is the following theorem.

\begin{theorem}\label{t9}
Let $\FF$ and $\FF'$ be two flags in $\F^n$ of length $r>2$.
\begin{enumerate}
\item If $\F$ is infinite and
$\sig(\FF)= (k, 1, l)$, $k,l>1$, then the semigroups
$S=\phi(\FF)$ and $T=\phi(\FF')$ are isomorphic
if and only if $\sig(\FF') = (k', 1, l')$, $k',l'>1$.
\item If $\F$ is infinite and $\sig(\FF)$ is different from
$(k, 1, l)$, $k,l>1$, then
the semigroups $S=\phi(\FF)$ and $T=\phi(\FF')$ are isomorphic
if and only if $\sig(\FF)=\sig(\FF')$.
\item If $\F$ is finite, then the semigroups $S=\phi(\FF)$ and
$T=\phi(\FF')$ are isomorphic if and only if
$\sig(\FF)=\sig(\FF')$.
\end{enumerate}
\end{theorem}

To prove this theorem we will need more notation and several
lemmas.

Let $\F$ be a field of an infinite cardinality $\gamma$ and
$V$  be a vector-space over $\F$ with a fixed basis,
$\{e_i|i=1,\dots,n\}$. For $v=\sum_{i=1}^{n}a_i e_i$ we set $\hat{v}=
(a_1,\dots,a_n)$ and denote by $\hat{V}$ the set of all $0\neq v\in V$,
such that $a_i=1$ for the minimal $i$ for which $a_i\neq 0$. For
$v\in\hat{V}$ set $l_v=\{\lambda v|\lambda\in\F^*\}$ and we have that
$V$ is a disjoint union of $\{0\}$ and all $l_v$, $v\in\hat{V}$.
Because of our choice of $\gamma$ we have $|l_v|=|V|=\gamma$.
Moreover, $|\hat{V}|=\gamma$ if $n>1$.

\begin{lemma}\label{ll1}
Let $M=(m_{i,j})_{i=1,\dots,k}^{j=1,\dots,l}$ be a $k\times l$
matrix of rank $1$ over $\F$. Then there exist
unique $\lambda\in\F^*$, $v\in\hat{\F^k}$ and
$w\in\hat{\F^l}$ such that $M=\lambda vw^t$.
\end{lemma}

\begin{proof}
Let $v'$ (resp. $w'$) be a non-zero column (resp. row) of $M$,
which exists, since $M$ has rank $1$. Then $v'=\lambda_1 v$ and
$w'=\lambda_2 w$ for some $\lambda_{1},\lambda_2\in\F$,
$v\in\hat{\F^k}$
and $w\in\hat{\F^l}$. Assume that $i'$ (resp. $j'$) is the number of
the first non-zero coordinate of $v$ (resp. $w$). Then, taking
$\lambda=m_{i',j'}$, we obviously get $M=\lambda vw^t$.
\end{proof}

\begin{proposition}\label{pp1}
Assume that $\F$ is infinite. Let $S$ and $T$ be some $3$-maximal
nilpotent semigroups of signature $(k,1,l)$
and $(m,1,n)$ with $k,l,m,n>1$ respectively. Then $S\simeq T$.
\end{proposition}

\begin{proof}
Let $\FF$ be the flag of $\F^{k+l+1}$, stabilized by $S$ and
$\FF'$ be the flag of $\F^{m+n+1}$, stabilized by $T$. Choose
also an $\FF$-basis, $\{e_i\}$, in $\F^{k+l+1}$, and an
$\FF'$-basis, $\{f_j\}$, in $\F^{m+n+1}$. Then the
elements of $S$ are matrices of the form
\begin{displaymath}
\left(\begin{array}{ccc}
0 & v & X \\
0 & 0 & w \\
0 & 0 & 0
\end{array}
\right),
\end{displaymath}
where $v\in\F^{k}$, $w^t\in\F^{l}$ and $X$ is a $k\times l$ matrix
over $\F$. Analogously, the elements of $T$ are matrices of the form
\begin{displaymath}
\left(\begin{array}{ccc}
0 & a & Y \\
0 & 0 & b \\
0 & 0 & 0
\end{array}
\right),
\end{displaymath}
where $a\in\F^{m}$, $b^t\in\F^{n}$ and $Y$ is an $m\times n$ matrix
over $\F$. We construct a bijection, $\phi:S\to T$, in several steps.

First, we construct bijections $\psi_1:\F^{k}\to \F^{m}$ and
$\psi_2:\F^{l}\to\F^{n}$ as follows: both $\psi_1$ and $\psi_2$ send
$0$ to $0$; then, as $k,l,m,n>1$, the sets $\hat{\F^{i}}$, $i=k,l,m,n$,
have the same cardinality and we can consider some bijections
$\psi'_1:\hat{\F^{k}}\to\hat{\F^{m}}$ and
$\psi'_2:\hat{\F^{l}}\to\hat{\F^{m}}$; the latter uniquely extend to
$\psi_1$ and $\psi_2$ by the properties $\psi_1(\lambda
x)=\lambda\psi_1(x)$, $\lambda\in\F$, $x\in \F^{k}$, and $\psi_2(\lambda
x)=\lambda\psi_2(x)$, $\lambda\in\F$, $x\in \F^{l}$.

Second, we decompose the set of all $k\times l$ (resp. $m\times n$)
matrices over $\F$ into a disjoint union of three subsets:
$\{0\}$, $M(1)$ and $M(2)$ (resp. $\{0\}$, $M(1)'$ and $M(2)'$), where
the first set contains the zero matrix, the second one contains all
matrices of rank $1$ and the last one
contains all other matrices. Since $\F$ is infinite,
we have $|M(1)|=|M(2)|=|M(1)'|=|M(2)'|$.
Let $\psi_3$ denote some bijection from $M(2)$ to $M(2)'$. The last
bijection $\psi_4:M(1)\to M(1)'$ is defined as follows: we take
$M\in M(1)$ and use Lemma~\ref{ll1} to write $M=\lambda v w^t$ for
some $\lambda\in\F$, $v\in\hat{\F^{k}}$ and $w\in\hat{\F^{l}}$;
then we set $\psi_4(M)=\lambda \psi_1(v)\psi_2(w)^t$. Lemma~\ref{ll1}
guarantees that $\psi_4$, as defined above, is a bijection. Let
$\psi_5$ be the bijection between the sets of all $k\times l$ and
$m\times n$ matrices, composed from $\psi_3$ and $\psi_4$.

Now we define $\phi:S\to T$ as follows:
\begin{displaymath}
S\ni\left(\begin{array}{ccc}
0 & v & X \\
0 & 0 & w \\
0 & 0 & 0
\end{array}
\right)\overset{\phi}{\mapsto}
\left(\begin{array}{ccc}
0 & \psi_1(v) & \psi_5(X) \\
0 & 0 & \psi_2(w^t)^t \\
0 & 0 & 0
\end{array}
\right)\in T.
\end{displaymath}

Finally, we are going to prove that $\phi$ is a homomorphism. In $S$
we have
\begin{displaymath}
\left(\begin{array}{ccc}
0 & v & X \\
0 & 0 & w \\
0 & 0 & 0
\end{array}
\right)
\left(\begin{array}{ccc}
0 & v' & X' \\
0 & 0 & w' \\
0 & 0 & 0
\end{array}
\right)=
\left(\begin{array}{ccc}
0 & 0 & vw' \\
0 & 0 & 0 \\
0 & 0 & 0
\end{array}
\right).
\end{displaymath}
Applying $\phi$ we get
\begin{displaymath}
\left(\begin{array}{ccc}
0 & \psi_1(v) & \psi_5(X) \\
0 & 0 & \psi_2(w^t)^t \\
0 & 0 & 0
\end{array}
\right)
\left(\begin{array}{ccc}
0 & \psi_1(v') & \psi_5(X') \\
0 & 0 & \psi_2({w'}^t)^t \\
0 & 0 & 0
\end{array}
\right)=
\left(\begin{array}{ccc}
0 & 0 & \psi_5(vw') \\
0 & 0 & 0 \\
0 & 0 & 0
\end{array}
\right),
\end{displaymath}
and the left-hand side equals
\begin{displaymath}
\left(\begin{array}{ccc}
0 & 0 & \psi_1(v)\psi_2({w'}^t)^t \\
0 & 0 & 0 \\
0 & 0 & 0
\end{array}
\right).
\end{displaymath}
Hence, we have to prove the equality $\psi_5(vw')=
\psi_1(v)\psi_2({w'}^t)^t$. This equality is an equality on
$\{0\}\cup M(1)'$ and follows directly from the definition of
$\psi_4$ (the $M(1)$-part of the map $\psi_5$). This completes
the proof.
\end{proof}

\begin{proposition}\label{ness1}
Let $T_1,$ $T_2$ be two isomorphic $r$-maximal nilpotent subsemigroups
of the semigroup $M(n,\F)$, $r \geq 3$,
$\sig(T_1)=(i_1, i_2, \dots, i_r)$ and
$\sig(T_2) = (j_1, j_2, \dots, j_r)$. Then $i_t = j_t$ for all
$t=2,\dots,r-1$.
\end{proposition}

To prove this statement we will need some notation. Recall that a
reflexive and transitive relation, $<$, on a set, X, is called a
preorder. A preorder, $<$, on $X$ induces the equivalence $\sim$
on $X$ defined by: $a \sim b$ if and only if  $a<b$ and $b<a$.
Then the relation $<$ becomes a partial order on the quotient set
$X/\sim$. By the {\em height} (resp. {\em depth}) of the preorder
$<$ we will understand the height (resp. depth) of the induced
partial order on $X/\sim$.

Let $T=\phi(\FF)$ be an $r$-maximal subsemigroup in $M(n,\F)$
(here $\FF$ is a flag of length $r$). For $A,B\in T$ we set $A\prec
B$ provided that $AC=0$ implies $BC=0$ for all $C\in T$, and we set
$A\ll B$ provided that $CA=0$ implies $CB=0$ for all $C\in T$. It is
clear that both $\prec$ and $\ll$ are preorders.

\begin{lemma}\label{lnnnnew1}
\begin{enumerate}
\item $A\prec B$ if and only if $(V_{r-1}\cap \ker(A))\subset
(V_{r-1}\cap \ker(B))$.
\item The set $M_0^\prec$ of all maximal elements with respect to
$\prec$ consists of all elements $A\in T$ such that $\ker(A)\supset
V_{r-1}$.
\item The set $M_i^\prec$ of all elements of depth $i$, $i=1,2$,
with respect to $\prec$ consists of all elements $A$ such that
$\dim(\ker(A)\cap V_{r-1})=\dim(V_{r-1})-i$.
\end{enumerate}
\end{lemma}

\begin{proof}
Let $A\prec B$ and suppose that there exists $v\in (V_{r-1}\cap
\ker(A))$ such that $v\not \in (V_{r-1}\cap \ker(B))$. Since $v\in
V_{r-1}$, there exists $C\in T$ such that $\langle
v\rangle=\imm(C)$. This guarantees $AC=0$ and $BC\neq 0$, a
contradiction. That $A\prec B$ follows from
$(V_{r-1}\cap \ker(A))\subset (V_{r-1}\cap \ker(B))$ is obvious.
This proves the first statement and the second statement follows
directly from the first one and Lemma~\ref{llnew}.

Let $A\in M_1^\prec$. Then $\ker(A)\not\supset V_{r-1}$ by the second
statement. Assume that $\dim(\ker(A)\cap V_{r-1})<\dim(V_{r-1})-1$ and
take arbitrary $0\neq v\in V_{r-1}\setminus \ker(A)$. Choose an
$\FF$-basis, $\{v_i\}$, of $\F^n$, which contains $v$ and a basis of
$\ker(A)\cap V_{r-1}$, and consider the linear operator $B$ defined
as follows: $B(v)=0$, $B(v_i)=A(v_i)$, $v_i\neq v$. It follows
directly form the definition and the first statement that $B\in T$ and
that $A\prec B$, moreover, $\dim(V_{r-1}\cap \ker(B))=\dim(V_{r-1}\cap
\ker(A))+1$. Hence $\ker(B)\not\supset V_{r-1}$ as $\dim(\ker(A)\cap
V_{r-1})<\dim(V_{r-1})-1$, and thus $B\not\in M_0^\prec$ by the second
statement. This implies that $A\not\in M_1^\prec$, a contradiction.
That all $A$ such that $\dim(\ker(A)\cap V_{r-1})=\dim(V_{r-1})-1$
belong to $M_1^\prec$ is obvious. This completes the proof for
$M_1^\prec$ and the arguments for $M_2^\prec$ are similar.
\end{proof}

\begin{lemma}\label{lnnnnew2}
\begin{enumerate}
\item $A\ll B$ if and only if $((V\setminus V_1)\cap \imm(A))\supset
((V\setminus V_1)\cap \imm(B))$.
\item The set $M_0^\ll$ of all maximal elements with respect to
$\ll$ consists of all elements $A$ such that $\imm(B)\subset
V_{1}$.
\item The set $M_i^\ll$ of all elements of depth $i$, $i=1,2$,
with respect to $\ll$ consists of all elements $A$ such that
$\dim(\imm(A))=\dim(\imm(A)\cap V_{1})+i$.
\end{enumerate}
\end{lemma}

\begin{proof}
Analogous to that of Lemma~\ref{lnnnnew1} with corresponding changes
and consideration of images instead of kernels.
\end{proof}

The main difficulty in the study of nilpotent subsemigroups in
$M(n,\F)$ up to isomorphism is that usual ranks of matrices are not
preserved under isomorphisms. To improve the situation, we introduce
the notion of a {\em super rank} (denoted by $\sr$), which we define
in the following way.

Denote $K_{1,0}=M_1^\prec\cap M_0^\ll$ and $K_{0,1}=M_0^\prec\cap
M_1^\ll$. Consider the set $M_1^\prec\cap M_1^\ll$ and define the
subset
\begin{displaymath}
K_{1,1}=\{A\in M_1^\prec\cap M_1^\ll\,:\, TAT\neq 0\}.
\end{displaymath}

An element, $A\in M(n,\F)$, is called {\em indecomposable} provided
that $A$ can not be decomposed into a product of two elements from
$T$. For an indecomposable element, $A\in M(n,\F)$, we will write
$\sr(A)=1$ provided that $A\in K_{1,1}\cup K_{1,0}\cup K_{0,1}$.

Define $K_{2,0}=M_2^\prec\cap M_0^\ll$ and $K_{0,2}=M_0^\prec\cap
M_2^\ll$. Consider the set $M_2^\prec\cap M_1^\ll$ and define the
subset
\begin{displaymath}
K_{2,1}=\{A\in M_2^\prec\cap M_1^\ll\,:\, TAT\neq 0\}.
\end{displaymath}
Consider the set $M_1^\prec\cap M_2^\ll$ and define the
subset
\begin{displaymath}
K_{1,2}=\{A\in M_1^\prec\cap M_2^\ll\,:\, TAT\neq 0\}.
\end{displaymath}
Consider the set $M_2^\prec\cap M_2^\ll$ and define the
subset
\begin{displaymath}
K_{2,2}=\{A\in M_2^\prec\cap M_2^\ll\,:\, TAT\neq 0\text{ and }TAT\neq
TBT\text{ for all }B, \sr(B)=1\}.
\end{displaymath}

For an indecomposable element, $A\in M(n,\F)$, we will write
$\sr(A)=2$ provided that $A\in K_{2,0}\cup K_{0,2}\cup K_{2,1}\cup
K_{1,2}\cup K_{2,2}$.

\begin{lemma}\label{lnnnnew3}
Let $T_1$ and $T_2$ be two isomorphic semigroups of signatures
$(i_1,\dots,i_r)$ and $(j_1,\dots,j_r)$ respectively and
$\phi:T_1\to T_2$ be an isomorphism. Let $A\in T_1$ be an
indecomposable element of super rank $i$, $i=1,2$. Then $\phi(A)$
is an indecomposable element of super rank $i$, $i=1,2$,
respectively.
\end{lemma}

\begin{proof}
Follows from the fact that isomorphism preserves indecomposable
elements and the sets $K_{u,v}$, where $0\leq u,v\leq 2$.
\end{proof}

Recall that for a semigroup, $S$, the notation $S^1$ denotes either
$S$ provided that $S$ has the identity or the semigroup $S\cup\{1\}$
in the opposite case. Our main aim in the following lemma is to
describe the super ranks for indecomposable elements in terms of the
multiplication.

\begin{lemma}\label{lnnnnew4}
Let $A\in T$ be an indecomposable element and $T^1AT^1\neq \{A,0\}$.
\begin{enumerate}
\item $\sr(A)=1$ if and only if there exists an element, $B\in T$,
of (usual) rank $1$, such that $CAD=CBD$ for all $C,D\in T^1$,
where either $C\neq 1$ or $D\neq 1$.
\item $\sr(A)=2$ if and only if $\sr(A)\neq 1$ and there exists an
element, $B\in T$, of (usual) rank $2$, such that $CAD=CBD$ for all
$C,D\in T^1$, where either $C\neq 1$ or $D\neq 1$.
\end{enumerate}
\end{lemma}

\begin{proof}
We prove the first statement and the second one is analogous. Let
$A=(a_{i,j})$ be an indecomposable element with $\sr(A)=1$. We start
with the necessity.

Let $\{v_i\}$ be an $\FF$-basis of $\F^n$.
If $A\in K_{1,0}$ then $M\in M_1^\prec$ and $M\in M_0^\ll$.
From Lemmas~\ref{lnnnnew1} and ~\ref{lnnnnew2} we get
$\dim(\ker(A)\cap V_{r-1})=\dim(V_{r-1})-1$ and
$\imm(A)\subset V_1$. Consider the element $B\in T$, defined as
follows: $B(v_i)=A(v_i)$, $i\leq \dim(V_{r-1})$, and $B(v_i)=0$ in
other cases. From $\dim(\ker(A)\cap V_{r-1})=\dim(V_{r-1})-1$ it
follows that $\rank(B)=1$ and it is obvious that $CAD=CBD$ for all
$C,D\in T^1$ such that either $C\neq 1$ or $D\neq 1$.

The case $A\in K_{0,1}$ can be handled by analogous argument.
Hence we consider the remaining case $A\in K_{1,1}$. From
Lemmas~\ref{lnnnnew1} and ~\ref{lnnnnew2} we get
$\dim(\ker(A)\cap V_{r-1})=\dim(V_{r-1})-1$ and
$\dim(\imm(A))=\dim(\imm(A)\cap V_1)+1$. Let $\{v_i\}$ be an
$\FF$-basis of $\F^n$. Consider the element
$B\in T$, defined via  $B(v_i)=A(v_i)$, $\dim(V_{1})< i\leq
\dim(V_{r-1})$, and $B(v_i)=0$ in other cases. Already the
condition $\dim(\ker(A)\cap V_{r-1})=\dim(V_{r-1})-1$
guarantees that $\rank(B)\leq 1$ and the condition $TAT\neq 0$ guarantees
$\rank(B)>0$. Therefore $\rank(B)=1$ and one easily gets
$CAD=CBD$ for all $C,D\in T^1$ such that either $C\neq
1$ or $D\neq 1$. This completes the proof of the necessity.

Let $A,B\in M(n,\F)$ be such that $T^1AT^1\neq \{A,0\}$,
$\rank(B)=1$ and $CAD=CBD$ for all $C,D\in T^1$ such that either
$C\neq 1$ or $D\neq 1$. Then we have
$\ker(A)\cap V_{r-1}=\ker(B)\cap V_{r-1}$ and $\imm(A)\setminus
(\imm(A)\cap V_1)=$ $\imm(B)\setminus
(\imm(B)\cap V_1)$. Hence $\dim(\ker(A)\cap V_{r-1})=$
$\dim(\ker(B)\cap V_{r-1})$. Since $\rank(B)=1$, we get that either
$\dim(\ker(A)\cap V_{r-1})=\dim(V_{r-1})$ and hence $A\in M_0^\prec$
or $\dim(\ker(A)\cap V_{r-1})=\dim(V_{r-1})-1$ and hence $A\in
M_1^\prec$. Analogously one gets $A\in M_0^\ll$ or $A\in M_1^\ll$. If
$A\in M_1^\prec\cap M_1^\ll$ then $T^1AT^1\neq \{A,0\}$ guarantees that
$A\in K_{1,1}$. This completes the proof.
\end{proof}

\begin{proof}[Proof of Proposition~\ref{ness1}]
Fix $1<s<r$. For $i=1,2$ consider the following sets:
\begin{gather*}
T_i(s)=\{A\in T_i\,:\, T_i^{s-2}AT_i^{r-s}\neq 0\text{ and }
\sr(A)=1\},\\
\tilde{T}_i(s)=T_i^{r-s}\cap\{B\in T_i\,:\, T^{s-1}B\neq 0\}.
\end{gather*}
Let $\phi:T_1\to T_2$ be an isomorphism. It follows from
Lemma~\ref{lnnnnew3} that $\phi(T_1(s))=T_2(s)$. Denote by $u(i)$,
$i=1,2$, the minimal positive integer for which there exists a subset,
$T'_i(s)\subset T_i(s)$, such that $|T'_i(s)|=u(i)$ and for all
$B\in \tilde{T}_i(s)$ there exists $C\in T'_i(s)$ such that $CB\neq
0$. It is clear that $u(1)=u(2)$ and to complete the proof it is
enough to show that $u(1)=i_s$.

First we show that $u(1)\geq i_s$. Indeed, let $H=\{A_1,\dots,A_j\}
\subset T_i(s)$ and $j< i_s$. Using Lemma~\ref{lnnnnew4}, we can find
the elements $B_1,\dots,B_j$ such that $\rank(B_l)=1$ for all $l$ and
$CA_lD=CB_lD$ for all $C,D\in T^1$, where either $C\neq
1$ or $D\neq 1$. But then the rank of arbitrary linear combination
of $B_l$'s does not exceed $j<i_s$. Since $\dim(V_s/V_{s-1})=i_s>j$ we
get that there exists $v\in V_s\setminus V_{s-1}$ such that $B_l(v)=0$
for all $l$. This means that $A_lC=B_lC=0$ for all $l$ if
$C\in T_i^{r-s}$ such that $\imm(C)=\langle v\rangle$. As
$T^{s-1}C=0$ and hence $C\in \tilde{T}_1(s)$, we get that $H$ does
not satisfy the necessary conditions.

Finally, one can construct a subset, $H\subset T_i(s)$, such that
$|H|=i_s$ and for all $B\in \tilde{T}_i(s)$ there exists $C\in H$
satisfying $CB\neq 0$, in the following way: let $0\neq v\in V_{s-1}$
and $\{v_1,\dots, v_{i_s}\}$ be arbitrary elements from $V_s$, which
are mapped onto a basis of $V_s/V_{s-1}$ under the canonical
projection. Extend $\{v_1,\dots, v_{i_s}\}\cup\{v\}$ to an
$\FF$-basis, $\mathcal{B}$ say, of $\F^n$. Let $A_l$, $l=1,\dots,i_s$,
be elements of rank $1$ in $M(n,\F)$, defined via
$A_l(v_l)=v$ and $A_l(w)=0$ for all $w\in
\mathcal{B}\setminus\{v_l\}$.  It is obvious that the set
$H=\{A_1,\dots, A_{i_s}\}\subset T$ satisfies the necessary
condition. This completes the proof.
\end{proof}

To proceed we will need the following lemma.

\begin{lemma} \label{posets}
Let $V$ and $V'$ be two vector spaces over a field. Let $V(2)$ and
$V'(2)$ denote the sets of all subspaces of $V$ and $V'$
respectively, whose dimensions do not exceed $2$. If $V(2)$ and
$V'(2)$ are isomorphic as posets with respect to inclusions then
$\dim V = \dim V'$.
\end{lemma}

\begin{proof}
It is enough to show that the poset of all subspaces  of $V$ (with respect to
inclusions) can be reconstructed from $V(2)$. To characterize 
$3$-dimensional subspaces we remark that every such a subspace is determined 
by some $2$-dimensional subspace, $\alpha$ say, and some $1$-dimensional
subspace, $l$ say, such that $l\not\subset \alpha$. In other words,
every pair $(\alpha,l)$, $l\not\subset \alpha$, defines a $3$-dimensional
subspace, which we will denote by $(\alpha,l)$, abusing notation. A
$1$-dimensional subspace, $m$, is contained in $(\alpha,l)$ if and only if
there exists a $1$-dimensional subspace, $k$, contained in $\alpha$, such that
$m$ is contained in the two-dimensional subspace, generated by $k$ and $l$.
Hence we can characterize all $1$-dimensional subspaces in $(\alpha,l)$ only 
in terms of elements from $V(2)$. A $2$-dimensional subspace $k$ is contained
in $(\alpha,l)$ if and only if all $1$-dimensional subsets of $k$ are contained
in $(\alpha,l)$. A subspace $k$ is contained in a subspace $m$ if and only if
all $1$-dimensional subsets of $k$ are contained in $m$. This implies that 
all subspaces of $V$ of dimension at most $3$ and their inclusions can 
be characterized only in terms of elements from  $V(2)$ and their 
inclusions. Hence the poset $V(2)$ completely determines the poset $V(3)$ or
all subspaces of $V$ of dimension at most $3$ with respect to inclusions. 
The proof is now easily completed by induction.
\end{proof}

Let $0\neq A\in T$ be a decomposable element. We will say that the
{\em super rank} $\sr(A)=1$ if there exists a decomposition,
$A=A_1A_2\dots A_l$, of $A$ into a product of indecomposable elements
and $j\in\{1,\dots,l\}$ such that $\sr(A_j)=1$. We will say that the
{\em super rank} $\sr(A)=2$ if $\sr(A)\neq 1$ and there exists a
decomposition, $A=A_1A_2\dots A_l$, of $A$ into a product of
indecomposable elements and $j\in\{1,\dots,l\}$ such that
$\sr(A_j)=2$.

\begin{lemma}\label{srdec}
Let $A\in T$ be decomposable.
\begin{enumerate}
\item $\rank(A)=1$ if and only if $\sr(A)=1$.
\item $\rank(A)=2$ if and only if $\sr(A)=2$.
\end{enumerate}
\end{lemma}

\begin{proof}
Let $\sr(A)=1$, $A=A_1A_2\dots A_l$ and $j\in\{1,\dots,l\}$ be such
that $\sr(A_j)=1$. Let $B$ be the element of rank $1$ such that
$CA_jD=CBD$ for all $C,D\in T^1$, where either $C\neq
1$ or $D\neq 1$. Then $A=A_1\dots A_{j-1}BA_{j+1}\dots A_l$ and hence
$\rank(A)\leq\rank(B)=1$. That $\rank(A)=1$ now follows from $A\neq
0$.

Conversely, let $\rank(A)=1$ and let $A=A_1A_2\dots A_l$ be a
decomposition of $A$ into indecomposable elements from $T$.
Let $v$ be a non-zero element in $\imm(A)$ and $\{v_i\}$ be an
$\FF$-basis of $\F^n$, such that $v\in\{v_i\}$. Consider the linear
operator $B$, defined as follows: $B(v)=v$ and $B(v_i)=0$, $v_i\neq
v$. Then  $A=A_1A_2\dots A_l=BA_1A_2\dots A_l$ and $\rank(B)=$
$\rank(BA_1)=1$. Moreover, from the construction of $B$ it follows
immediately that $BA_1\in T$. If $BA_1$ is indecomposable, we get
$\rank(BA_1)=1$ and $T^1 BA_1T^1\ni A$ and hence does not coincide
with $\{BA_1,0\}$. This implies $\sr(BA_1)=1$. Otherwise we can
substitute the element $A$ by the element $BA_1$ and apply the same
arguments. This procedure will certainly stop after not more than
$n$ steps, since $T$ is a nilpotent semigroup of nilpotency degree
$\leq n$ (see Corollary~\ref{cnew1}).

The proof of the second statement is analogous.
\end{proof}

\begin{proof}[Proof of Theorem~\ref{t9}]
Let $T_1$ and $T_2$ be two $r$-maximal nilpotent subsemigroup in
$M(n,\F)$. Assume that $\sig(T_1)=(i_1,\dots,i_r)$ and
$\FF:0=V_0\subset V_1\subset \dots$ be the corresponding flag in
$\F^n$. Assume that $\sig(T_2)=(j_1,\dots,j_r)$ and
$\FF:0=W_0\subset W_1\subset \dots$ be the corresponding flag in
$\F^n$. The strategy of our proof is to show that all $i_s$ can be
characterized in terms of some invariants, preserved under any
isomorphism. If $r>3$ and  $1<s<r$, this was done in
Proposition~\ref{ness1}. Hence here we have to consider only remaining
cases.

{\bf Case 1.} $i_1,i_r,\sum_{s=2}^{r-1}i_s\geq 2$. Then it is easy to
see that for every $V'\subset V_1$, $\dim(V')=2$, there exists a
decomposable element, $A$, of rank $2$, such that $\imm(A)=V'$. Let
$\phi:T_1\to T_2$ be an isomorphism. Denote by $\Ann(T_1)$ the
two-sided annihilator of $T_1$. First we remark that the inequality
$\sum_{s=2}^{r-1}i_s\geq 2$ implies that  for a fixed $i=1,2$ and for
$A,B\in \Ann(T_1)$, $\rank(A)=\rank(B)=i$, the condition
$\imm(A)=\imm(B)$ is equivalent to the existence of $C$,
$\rank(C)=i$, such that $A=CD_1$ and $B=CD_2$ for some
$D_1$ and $D_2$. In particular, for this $C$ we will have
$\imm(C)=\imm(A)=\imm(B)$. Hence, $\imm(A)=\imm(B)$ is equivalent to
$\imm(\phi(A))=\imm(\phi(B))$. Now from Lemmas~\ref{srdec} and
\ref{lnnnnew3} it follows that
$\dim(\imm(\phi(A)))=\dim(\imm(A))$.

Let $A,B$ be two elements from $\Ann(T_1)$, $\rank(A)=1$,
$\rank(B)=2$. Then it is easy to see that
$\imm(A)\subset \imm(B)$ is equivalent to the
condition that for arbitrary decomposition $B=CD_1$, with $C$
indecomposable and $\sr(C)=2$, there exists $D_2$ such that
$A=CD_2$. This implies that $\imm(A)\subset \imm(B)$ is equivalent to
$\imm(\phi(A))\subset \imm(\phi(B))$. Hence $\phi$ induces an
isomorphism from the partially ordered set of all subsets of dimension
$\leq 2$ in $V_1$ to the partially ordered set of all subsets of
dimension $\leq 2$ in $W_1$. Now Lemma~\ref{posets} guarantees
$\dim(V_1)=\dim(W_1)$.

{\bf Case 2.} $\sum_{s=2}^{r-1}i_s\geq 2$ and $i_1=1$. This situation
is obviously characterized by the fact that
$K_{2,0}=\varnothing$. Hence $i_1=1$ implies $j_1=1$. Moreover, going
to the opposite semigroup, we also get $i_r=1$ implies $j_r=1$ in the
case $\sum_{s=2}^{r-1}i_s\geq 2$.

{\bf Case 3.} $r=3$, $i_2=1$ and $i_1=1$. This case is characterized
by the fact that there exists an indecomposable element, $A$, such
that $\sr(A)=1$ and $AT$ contains all decomposable elements. Hence in
this case we get $j_1=1$. Analogously $i_3=1$ implies $j_3=1$.

{\bf Case 4.} $r=3$, $i_2=1$, $i_1,i_3>1$ and $\F$ is infinite. From
previous cases we get that $j_1,j_3>1$ as well. By
Proposition~\ref{pp1}, $T_2\simeq T_1$ for arbitrary $j_1,j_3>1$.

{\bf Case 5.} $r=3$, $i_2=1$, $i_1,i_3>1$ and $|\F|=q=p^l<\infty$.
One easily calculates that the number of decomposable elements in
$\Ann(T_1)$ ($i_1\times i_3$-matrices of rank $\leq 1$) equals
$q^{i_1+i_3-1}$. Hence, the number of indecomposable elements in
$\Ann(T_1)$ ($i_1\times i_3$-matrices of rank $>1$) equals
$q^{i_1+i_3}-q^{i_1+i_3-1}$. Hence $T_1\simeq T_2$ implies
$q^{i_1+i_3}-q^{i_1+i_3-1}=q^{j_1+j_3}-q^{j_1+j_3-1}$ and from
$i_1+i_3=j_1+j_3=n-1$ we get $i_1i_3=j_1j_3$. Further, one can easily
see that the number of elements in the right annihilator of $T_1$
equals $q^{i_1(i_2+i_3)}$. From $T_1\simeq T_2$ we get
$q^{i_1(i_2+i_3)}=$ $q^{j_1(j_2+j_3)}$ and, finally, from $i_2=j_2$
and from $i_1i_3=j_1j_3$ we derive $i_1=j_1$. This also implies
$i_3=j_3$ and completes the proof.
\end{proof}

Careful analysis of the proof above shows that, in fact, the following
statement is true.

\begin{corollary}\label{cccnnn}
Let $r>3$ and $T_1$, $T_2$ be some $r$-maximal nilpotent semigroups
(possibly in different $M(i,\F)$). Then $T_1\simeq T_2$ if and only if
$\sig(T_1)=\sig(T_2)$.
\end{corollary}

\section{Appendix: Isolated and completely isolated subsemigroups
(in the case of a finite $\F$)}\label{sn5}

For $i=0,\dots,n$ we denote by $D_i$ the set of all $A\in M(n,\F)$
such that the rank of $A$ equals $i$. The sets $D_i$, $i=0,\dots,n$,
constitute an exhaustive list of $\mathcal{D}$-classes in $M(n,\F)$,
see \cite[Chapter~2]{Ok2}. For $i=0,\dots,n$ we also denote by $I_i$
the two-sided ideal  $D_0\cup\dots\cup D_i$ of $M(n,\F)$.
We start with the following lemma, which we
will often use in the sequel.

\begin{lemma}\label{l10}
For every positive integer $k$, $1 \leq k \leq n-1$, the ideal
$I_k$ is generated as a semigroup by $D_k$.
\end{lemma}

\begin{proof}
Let $S=\langle D_k\rangle$. Since $I_k=\cup_{j\leq k} D_j$, by
inductive arguments it is enough to show that  $D_{k-1}\subset S$.
Let $A\in D_{k-1}$, $V_1=\imm(A)$, $V_2=\ker(A)$. Then
$\dim(\ker(A))\geq 2$. Let $v_1,\dots, v_i$ be a basis of
$\ker(A)$, $v_{i+1},\dots,v_n$ a basis of $\imm(A)$, and
$w_{i+1},\dots,w_n$ be arbitrary such that $v_{j}=A(w_j)$, $j>i$.
Note that this implies that $w_{i+1},\dots,w_n$ are linearly
independent and $v_1,\dots,v_i,w_{i+1},\dots, w_n$ form the basis
of $\F^n$. Choose further two arbitrary elements $w,u$, such that
$\dim(\langle w,u,V_1\rangle)= \dim(V_1)+2$. Note that
$i=n-k+2\geq 2$ by choice of $k$. Define the linear operator $B_1$
as follows: $B_1(v_l)=0$, $l=1,\dots,i-1$; $B_1(v_i)=w$;
$B_1(w_l)=v_l$, $l>i$. Let $B_2$ be a linear operator of rank $k$,
satisfying $B_2(u)=u$, $B_2(w)=0$, $B_2(v_l)=v_l$, $l>i$. Clearly
such $B_2$ exists and $B_2B_1=A$. This completes the proof.
\end{proof}

Recall that a subsemigroup, $T$, of a semigroups, $S$, is called {\em
completely isolated} provided that $xy\in T$ implies $x\in T$ or $y\in
T$ for arbitrary $x,y\in S$. Further, $T$ is called {\em isolated}
provided that $x^n\in T$ implies $x\in T$ for every $x\in S$.

\begin{proposition}\label{p1sn5}
Assume that $\F$ is finite. Then the only completely isolated
subsemigroups of $M(n,\F)$ are $M(n,\F)$, $\GL(n,\F)$ and $I_{n-1}=
M(n,\F)\setminus \GL(n,\F)$.
\end{proposition}

\begin{proof}
It is obvious that the subsemigroups $M(n,\F)$, $\GL(n,\F)$ and
$I_{n-1}$ are completely isolated. Let $S$ be a completely isolated
subsemigroup of $M(n,\F)$. The semigroup $S\cap \GL(n,\F)$ is finite
and hence must contain an idempotent provided that it is not
empty. Since $\GL(n,\F)$ contains only one idempotent we get that
either $S\supset \GL(n,\F)$ or $S\cap \GL(n,\F)=\varnothing$.

Assume that $S\cap I_{n-1}\neq\varnothing$. By Lemma~\ref{l10},
$I_{n-1}$ is generated by $D_{n-1}$. Hence every element from
$S\cap I_{n-1}$ can be written as a product of elements from
$D_{n-1}$. Since $S$ is completely isolated, we conclude that
$S\cap D_{n-1}$ is not empty. Let $A\in S\cap D_{n-1}$. For every
$B\in D_{n-1}$ there exist $C,D\in \GL(n,\F)$ such that $A=CBD$.
If $S\supset \GL(n,\F)$ then $B=C^{-1}AD^{-1}\in S$. Otherwise
$B\in S$ because $S$ is completely isolated. Thus $S\supset
D_{n-1}$. This means that $S\supset I_{n-1}=M(n,\F)\setminus
\GL(n,\F)$ by Lemma~\ref{l10} and the statement follows.
\end{proof}

Let $\mathcal{A}$ and $\mathcal{B}$ denote some non-empty sets of
$n-1$-dimensional and $1$-dimensional subspaces in $\F^n$
respectively. Assume that $V_1\not\subset V_2$ for every
$V_1\in\mathcal{B}$ and every $V_2\in\mathcal{A}$. Set
\begin{displaymath}
S(\mathcal{A},\mathcal{B})=\{A\in M(n,\F)\,:\, \imm(A)\in
\mathcal{A}, \ker(A)\in \mathcal{B}\}.
\end{displaymath}

\begin{lemma}\label{ll1sn5}
For every $\mathcal{A}$ and $\mathcal{B}$ as above the set
$S(\mathcal{A},\mathcal{B})$ is an isolated subsemigroup in
$M(n,\F)$.
\end{lemma}

\begin{proof}
Let $A,B\in S(\mathcal{A},\mathcal{B})$. Since $\imm(B)\cap
\ker(A)=0$, we get that $\ker(AB)=\ker(B)$ and
$\imm(AB)=\imm(A)$. This means that $AB\in S(\mathcal{A},\mathcal{B})$
and thus $S(\mathcal{A},\mathcal{B})$ is a subsemigroup in $M(n,\F)$.

Let $A\in M(n,\F)$ be such that $A^k\in S(\mathcal{A},\mathcal{B})$
for some positive integer $k$. Then $\rank(A^k)=n-1$ and hence
$\rank(A)=n-1$. This implies $\ker(A)=\ker(A^k)\in\mathcal{B}$ and
$\imm(A)=\imm(A^k)\in\mathcal{A}$. Therefore $A\in
S(\mathcal{A},\mathcal{B})$, which completes the proof.
\end{proof}

For any decomposition $\F^n=V_1\oplus V_2$, we denote by
$e(V_1,V_2)$ the uniquely defined idempotent from $M(n,\F)$, which
represents the linear operator of projection on $V_1$ along
$V_2$. Note that the set of all possible $e(V_1,V_2)$ constitute an
exhaustive list of idempotents in $M(n,\F)$.

\begin{lemma}\label{ll2sn5}
Let $\F$ be finite and $S$ be an isolated semigroup in $M(n,\F)$.
Assume that one of the following conditions is satisfied:
\begin{enumerate}
\item $S$ contains $0$.
\item $S$ contains two idempotents $e(V_1,V_2),e(V'_1,V'_2)$
of rank $n-1$, such that $V_2\subset V'_1$.
\item $S$ contains some idempotent of rank $\leq n-2$.
\end{enumerate}
Then $S\supset I_{n-1}$.
\end{lemma}

\begin{proof}
First we consider the case when $S$ contains $0$.
Since $S$ is isolated and contains $0$, it must contain all nilpotent
matrices. Let $V_1$ and  $V_2$ be two non-trivial
subspaces in $\F^n$ such that
$\F^n=V_1\oplus V_2$. Choose a basis, $v_1,\dots,v_i$, of $V_1$ and a
basis, $v_{i+1},\dots,v_n$, of $V_2$, and consider the linear operators
$A$ and $B$ defined as follows: $A(v_j)=v_{j-1}$, $j=2,\dots,i+1$,
$A(v_j)=0$ otherwise; $B(v_j)=v_{j+1}$, $j=1,\dots,i$, $B(v_j)=0$
otherwise. Clearly, both $A$ and $B$ are nilpotent and hence are
contained in $S$. Thus $S\ni AB=e(V_1,V_2)$, which means that $S$
contains all non-invertible idempotents. Since for a finite $\F$
some power of every element in $I_{n-1}$ is an idempotent, we 
finally get that $S\supset I_{n-1}$.

Now we consider the case when $S$ contains two idempotents 
$e(V_1,V_2),e(V'_1,V'_2)$ of rank $n-1$, such that $V_2\subset V'_1$.
Because of the first part of the lemma proved above it is enough to show 
that $S$ contains $0$. Let $v_1,\dots,v_{n-1}$ be a basis of $V'_1$, 
such that $v_{n-1}$ is a basis of $V_2$, and let $v_n$ be a basis of 
$V'_2$. Consider the linear operator $A$ defined as follows: 
$A(v_i)=v_{i+1}$, $i<n-1$, $A(v_{n-1})=v_1$ and $A(v_n)=0$. Then 
$A^{n-1}=e(V'_1,V'_2)\in S$ and
hence $A\in S$, since $S$ is isolated. Therefore $B_i=A^i e(V_1,V_2)
A^{n-1-i}\in S$ for every $i=1,\dots,n-1$. But $B_i(v_j)=v_j$, $j\neq
n,i$, and $B_i(v_n)=B_{i}(v_i)=0$. Hence $0=B_1B_2\dots B_{n-1}\in S$.

Finally, we consider the case when $S$ contains some 
idempotent of rank $\leq n-2$. Because of the second part of the lemma 
proved above it is enough to show that $S$ contains two
idempotents, say $e,f$, of rank $n-1$, such that $\ker(e)\in \imm(f)$.
Let $e(V_1,V_2)\in S$ be a non-invertible idempotent and
$\dim(V_2)\geq 2$. Choose a basis, $v_1,\dots,v_i$, of $V_1$ and a
basis $v_{i+1},\dots,v_n$, of $V_2$, and consider the linear
operators $A$ and $B$ defined as follows: $A(v_j)=v_j$, $j\leq i$,
$A(v_{i+1})=0$, $A(v_j)=v_{j-1}$, $j>i+1$; $B(v_j)=v_j$, $j\leq i$,
$B(v_{n})=0$, $B(v_j)=v_{j+1}$, $j=i+1,\dots,n-1$. Then
$A^{\dim(V_2)}=B^{\dim(V_2)}=e(V_1,V_2)\in S$ and thus $A,B\in S$ as
$S$ is isolated. But this means that $AB,BA\in S$. But both $AB$ and
$BA$ are idempotents of rank $n-1$ and by construction $\ker(AB)=
\langle v_n\rangle\subset \imm(BA)$. This completes the proof.
\end{proof}

\begin{theorem}\label{t2sn5}
Assume that $\F$ is finite. Then the only isolated subsemigroups
of the semigroup $M(n,\F)$ are $M(n,\F)$, $\GL(n,\F)$, $I_{n-1}$
and all semigroups $S(\mathcal{A},\mathcal{B})$, where
$\mathcal{A}$ and $\mathcal{B}$ are some non-empty sets of
$n-1$-dimensional and $1$-dimensional subspaces in $\F^n$
respectively, such that $V_1\not\subset V_2$ for every
$V_1\in\mathcal{B}$ and every $V_2\in\mathcal{A}$.
\end{theorem}

\begin{proof}
Since every completely isolated subsemigroup is isolated, it follows
from Lemma~\ref{ll1sn5} that all the semigroups listed above are indeed
isolated. Let now $S$ be an isolated semigroup of $M(n,\F)$.
By the same arguments as in Proposition~\ref{p1sn5}, we get that
either $S\supset \GL(n,\F)$ or $S\cap \GL(n,\F)=\varnothing$. We
consider these two cases separately.

{\bf Case 1.} $S\supset \GL(n,\F)$. If $S=\GL(n,\F)$, we are done.
So, we can assume that $S\cap
I_{n-1}\neq\varnothing$. Then $S\cap D_k\neq\varnothing$ for some $k$
and hence $S\supset D_k$ since $S\supset \GL(n,\F)$. Lemma~\ref{l10}
now implies $S\supset I_k$, in particular, $S\ni 0$. Thus
$S\supset I_{n-1}$ by Lemma~\ref{ll2sn5} and we get $S=M(n,\F)$ in
this case.

{\bf Case 2.} $S\cap \GL(n,\F)=\varnothing$. Then $S$ contains a
non-invertible idempotent. If $S$ contains an idempotent of rank $\leq
n-2$, then $S\supset I_{n-1}$ according to Lemma~\ref{ll2sn5} and thus
$S=I_{n-1}$. Hence we can now assume that all idempotents of $S$ have
rank $n-1$. Let $E(S)=\{e_1,\dots,e_m\}$ and set
$\mathcal{A}=\{\imm(e_1),\dots,\imm(e_m)\}$,
$\mathcal{B}=\{\ker(e_1),\dots,\ker(e_m)\}$. From Lemma~\ref{ll2sn5}
it follows that $\ker(e_j)\not\subset \imm(e_{j'})$ for all
$j,j'$. Let $V_1\in \mathcal{A}$ and $e_{\alpha}\in E(S)$ be such that
$V_1=\imm(e_\alpha)$. Let $V_2\in \mathcal{B}$ and $e_{\beta}\in E(S)$
be such that $V_2=\ker(e_\beta)$. Then $\imm((e_{\alpha}e_{\beta})^i)=V_1$
and $\ker((e_{\alpha}e_{\beta})^i)=V_2$ for all $i$, in particular,
for that $i$, for which the element $(e_{\alpha}e_{\beta})^i$ belongs
to $E(S)$. This proves that $E(S)\ni e(V_1,V_2)$ for all $V_1\in
\mathcal{A}$ and all $V_2\in \mathcal{B}$. Take now arbitrary $A\in
S(\mathcal{A},\mathcal{B})$. Since
$A^i\in E(S(\mathcal{A},\mathcal{B}))=E(S)$ for some $i$, we get
that $A\in S$ and hence $S\supset S(\mathcal{A},\mathcal{B})$. On
the other hand, the fact that $S(\mathcal{A},\mathcal{B})$ is isolated
says that $A\not\in S(\mathcal{A},\mathcal{B})$ implies that the
corresponding idempotent $A^i$ can not belong to
$E(S)=E(S(\mathcal{A},\mathcal{B}))$ since either
$\imm(A^i)\not\in\mathcal{A}$ or $\ker(A^i)\not\in\mathcal{B}$. From this
we derive $S=S(\mathcal{A},\mathcal{B})$ and complete the proof.
\end{proof}


\begin{center}
\bf Acknowledgments
\end{center}

This paper was essentially written during the visits of the first
author to Uppsala University, which were supported by The Royal Swedish
Academy of Sciences and The Swedish Institute. The financial support
of The Academy. The Swedish Institute and the  hospitality of
Uppsala University are gratefully acknowledged.
We also would like to thank Prof. O.~Ganyushkin for useful
conversations. We are in debt to Prof. T.~Ekholm and Prof. O.~Viro for
their help regarding the proof of  Lemma~\ref{posets}. We thank
Prof. J.~Okni\'nski and Prof. O.~Artemovych for their comments on the
first version of the paper. We are especially in debt to the referee for
many very valuable comments.

\vspace{1cm}

\noindent
V.M.: Department of Mathematics, Uppsala University, SE 471 06,
Uppsala, SWEDEN, e-mail: {\em mazor\symbol{64}math.uu.se}
\vspace{0.5cm}

\noindent
G.K.: Department of Mechanics and Mathematics, Kyiv Taras Shevchenko
University, 64, Volodymyrska st., 01033, Kyiv, UKRAINE,\\
e-mail: {\em akudr\symbol{64}univ.kiev.ua}

\end{document}